\documentclass[11pt,a4paper,twoside,leqno]{article}

\usepackage{times}
\usepackage{amsmath,amssymb,amsthm,amsfonts,epsfig}
\usepackage{a4wide}
\usepackage{exscale}   
\usepackage{enumerate}
\usepackage{latexsym}

\usepackage[all]{xy} 

\usepackage{stmaryrd}

\newcommand{\restr}{\rvert}   
\newcommand{\abs}[1]{\left\lvert#1\right\rvert}   
\newcommand{\norm}[1]{\left\lVert#1\right\rVert}

\newcommand{\Cont}{{\mathcal C}} 
\DeclareMathOperator{\supp}{supp}
\DeclareMathOperator{\top2}{top}

\DeclareMathOperator{\Lin}{L} 
\DeclareMathOperator{\Leb}{L} 

\DeclareMathOperator{\Ind}{Ind} 

\DeclareMathOperator{\Komp}{K} 

\DeclareMathOperator{\KTh}{K}
\DeclareMathOperator{\KK}{KK}

\newcommand{\C}{\ensuremath{{\mathbb C}}}
\newcommand{\E}{\ensuremath{{\mathbb E}}}

\newcommand{\R}{\ensuremath{{\mathbb R}}}

\newcommand{\mA}{\ensuremath{{\mathcal A}}}

\newcommand{\mG}{\ensuremath{{\mathcal G}}}

\theoremstyle{plain}

\newtheorem {theorem} {Theorem}[section]
\newtheorem {lemma}[theorem] {Lemma}
\newtheorem {proposition} [theorem]{Proposition}
\newtheorem {corollary} [theorem]{Corollary}

\newtheorem* {theorem*} {Theorem}
\newtheorem* {proposition*} {Proposition}
\newtheorem* {lemma*}{Lemma}

\theoremstyle{definition}

\newtheorem {definition} [theorem]{Definition}
\newtheorem {defprop} [theorem]{Definition and Proposition}

\newtheorem {remark} [theorem]{Remark}

\DeclareMathOperator{\KKban}{KK^{\ban}}

\DeclareMathOperator{\ban}{ban}

\newcommand{\rmd}{{\, \mathrm{d}}}

\newcommand{\ketbra}[2]{\big|#1\big\rangle\big\langle #2 \big|}

\newcommand{\gG}{\mG} 

\newcommand{\EbanW}[1]{\E^{\ban}_{#1}} 
\newcommand{\KKbanW}[1]{\KK^{\ban}_{#1}} 

\newcommand{\ContSect}{\Gamma}

\begin{document}

\title{The Bost conjecture, open subgroups and groups acting on trees}
\author{Walther Paravicini}
\date{May 31, 2008}
\maketitle

\begin{abstract}
\noindent The Bost conjecture with C$^*$-algebra coefficients for locally compact Hausdorff groups passes to open subgroups. We also prove that if a locally compact Hausdorff group acts on a tree, then the Bost conjecture with C$^*$-coefficients is true for the group if and only if it is true for the stabilisers of the vertices.

\medskip


\medskip

\noindent \textsc{AMS 2000} \textit{Mathematics subject classification:} Primary 19K35; 46L80; 22A22; 43A20
\end{abstract}

\noindent If $H$ is a closed subgroup of a locally compact Hausdorff group $G$, then the group $H$ is equivalent, as a groupoid, to the transformation groupoid $\gG:=G \ltimes G/H$. Given an $H$-C$^*$-algebra $B$ one can construct a canonical induced $G$-C$^*$-algebra $\Ind_H^G B$, but using the theory of equivalences of groupoids one can also construct a $\gG$-C$^*$-algebra $\Ind_H^{\gG} B$. The difference between the two is that $\Ind_H^{\gG} B$ keeps track of an additional $G/H$-fibration of the induced algebra, i.e., a compatible $\Cont_0(G/H)$-module structure.

A possible blueprint for a proof of the (already known) fact that the Baum-Connes conjecture passes to closed subgroups is the following: Firstly, prove that the Baum-Connes conjecture is stable under the passage to equivalent groupoids, e.g. from $H$ to $G\ltimes G/H$; secondly, analyse to what extend one can ``forget'' additional fibrations, e.g. the above-mentioned $\Cont_0(G/H)$-structure on $\Ind_H^{\gG} B$. It turns out that the following diagram is commutative
\[
\xymatrix{
\KTh^{\top2}_*(H, B) \ar[d]_{\cong} \ar[r] & \KTh_*( B \rtimes_r H) \ar[d]^{\cong}\\
\KTh^{\top2}_*(\gG , \Ind_H^{\gG}B) \ar[d]_{\cong} \ar[r] & \KTh_*( \Ind_H^{\gG} B \rtimes_r \gG) \ar[d]^{\cong}\\
\KTh^{\top2}_*(G, \Ind_H^G B)  \ar[r] & \KTh_*( \Ind_H^G B \rtimes_r G)
}
\]
The horizontal arrows in this diagram are the respective Baum-Connes assembly maps. The vertical arrow on the upper left-hand side is given by an induction homomorphism, the vertical arrow on the upper right-hand side stems from a Morita equivalence. The lower left vertical arrow is the so-called ``forgetful map'' and the remaining arrow is an isomorphism in $\KTh$-theory that corresponds to an isomorphism of the involved C$^*$-algebras.

Whereas most of the isomorphisms appearing in this diagram emerge quite naturally from the theory, the fact that the lower left homomorphism is an isomorphism is non-trivial (see \cite{CEO:03}, Theorem~0.1, or \cite{ChaEch:01:Permanence}, Theorem~2.2). The diagram shows that the Baum-Connes conjecture with coefficients for $G$ implies the Baum-Connes conjecture with coefficients for $H$.

In the present work, the same fundamental plan of attack is used to show that also the $\Leb^1$-version of the Baum-Connes conjecture (with C$^*$-algebra coefficients), called the Bost conjecture, passes to subgroups, at least if they are \emph{open}. The relevant diagram is the following, where the horizontal arrows are the Bost assembly maps constructed in \cite{Lafforgue:06}:
\[
\xymatrix{
\KTh^{\top2}_*(H, B) \ar[d]_{\cong} \ar[r] & \KTh_*( \Leb^1(H, B) ) \ar[d]^{\cong}\\
\KTh^{\top2}_*(\gG , \Ind_H^{\gG}B) \ar[d]_{\cong} \ar[r] & \KTh_*( \Leb^1(\gG, \Ind_H^{\gG} B ))  \\
\KTh^{\top2}_*(G, \Ind_H^G B)  \ar[r] & \KTh_*( \Leb^1(G, \Ind_H^G B )) \ar[u]_{\psi_*}
}
\]
The upper square in this diagram is obtained as a special case of the induction theorem for equivalent groupoids proved in \cite{Paravicini:07:Induction:arxiv}. The vertical isomorphism in the lower left-hand corner is the same as in the first diagram. Note that the reason why we are not able to treat Banach algebra coefficients instead of C$^*$-coefficients in the present paper is that we do not know whether the lower left-hand arrow is an isomorphism for Banach algebra coefficients (and $\KTh^{\top2,\ban}$ instead of $\KTh^{\top2}$). On the other hand, we do know that the vertical arrow labelled with $\psi_*$ is an isomorphism if $H$ is open or compact, and this even holds for arbitrary Banach algebra coefficients. Therefore, the Bost conjecture passes to open subgroups. If one could prove that $\psi_*$ is an isomorphism for general closed subgroups $H$, then the Bost conjecture would also pass to closed subgroups.

\medskip

Section~\ref{Section:TheForgetfulMap} introduces the pushforward of fields of Banach spaces, by which we implement the above-mentioned process of forgetting additional fibrations over locally compact spaces. We show that this forgetful map lifts to Banach algebras and also to $\KKban$-cycles, so we obtain a forgetful homomorphism on the level of $\KKban$-theory. We show that it is compatible with the descent construction and to some extend also with the Bost assembly map. All of this is done for arbitrary unconditional completions instead of $\Leb^1$, and for arbitrary actions of groupoids on arbitrary spaces instead of group actions on homogeneous spaces.

An analysis of the $\Leb^1$-algebras that arise in the study of group actions on locally compact spaces reveals in Section~\ref{Section:TheSpecialCaseOfL1} technical subtleties which so far prevent us from proving that the homomorphism $\psi_*$ in the above diagram is an isomorphism in general. However, in Section~\ref{Section:Subgroups}, we carry out the program sketched above and show that the $\Leb^1$-version of the Bost conjecture for a second countable locally compact Hausdorff groups passes to open subgroups.

In Section~\ref{Section:Tree}, we show that if a second countable locally compact Hausdorff groups acts on a tree, then the $\Leb^1$-version Bost conjecture with C$^*$-coefficients is true for the group if and only if it is true for the stabilisers of the vertices.

It is perhaps worth pointing out that at certain places we confined ourselves to treat just the $\Leb^1$-version of the Bost conjecture but it should be possible to extend the methods implemented for $\Leb^1$ to one or the other more general unconditional completion.

\medskip

I would like to thank Siegfried Echterhoff for his constructive and encouraging advice and for suggesting the inheritance properties of the Bost conjecture as a research topic; although the conclusions of the present article are newer, the structural results that are necessary to obtain them are indeed already contained in the doctoral thesis \cite{Paravicini:07}. Also, I would like thank Hervé Oyono-Oyono who has explained me how my results on induction should admit to treat the case of group actions on a tree.

This research has been supported by the Deutsche Forschungsgemeinschaft (SFB 478).

\smallskip

Notation: All Banach spaces and Banach algebras that appear in this article are supposed to be complex. As far as groupoids, $\KKban$-theory for Banach algebras and groupoids etc.\ are concerned we refer to the definitions in \cite{Lafforgue:06} and \cite{Paravicini:07:Induction:arxiv}.

\section{The forgetful map in the Banach algebra context}\label{Section:TheForgetfulMap}

In \cite{CEO:03}, the following results are shown: if $G$ is second countable locally compact Hausdorff group, $Y$ is a second countable locally compact Hausdorff $G$-space and $B$ is a separable $G\ltimes Y$-C$^*$-algebra, then there is an isomorphism $B \rtimes_r (G \ltimes Y) \cong \ContSect_0(Y,B) \rtimes_r G$ and a commutative diagram of the following form:
\begin{equation}\label{Diagram:CEO:BaumConnes}
\xymatrix{
\KTh^{\top2}_*(G\ltimes Y; B) \ar[rr] \ar[d]_{\cong} && \KTh_*(B \rtimes_r (G \ltimes Y)) \ar[d]^{\cong}\\
\KTh^{\top2}_*(G; \ContSect_0(Y,B)) \ar[rr] && \KTh_*(\ContSect_0(Y,B) \rtimes_r G)
}
\end{equation}
where the horizontal arrows are the respective Baum-Connes assembly maps. Note that it is a non-trivial fact that the vertical arrow on the left-hand side, called the \emph{forgetful map}, is an isomorphism (actually, this is the main result of \cite{CEO:03}).

In this section, we construct a version of this forgetful map for groupoids and Banach algebras, without being able to prove that it is an isomorphism. However we are able to prove that the Banach algebraic version (with unconditional completions) of the above diagram is commutative. From this it is easy to switch over to C$^*$-algebra coefficients and unconditional completions, where we have an isomorphism on the left-hand side. That is why the main applications of this article will require C$^*$-coefficients. It would probably not be much easier to consider C$^*$-algebra coefficients right from the start, and our approach prepares the ground for further studies of the Bost conjecture with Banach algebra coefficients.

For the notation concerning fields of Banach spaces and groupoids we refer the reader to \cite{Lafforgue:06} and \cite{Paravicini:07:Induction:arxiv}.

\subsection{The forgetful map for equivariant fields of Banach spaces etc.} \label{Subsection:ThePushforwardConstruction}

If $X$ and $Y$ are locally compact Hausdorff spaces and $p\colon Y\to X$ is continuous, then we want to know how to transform a field of Banach spaces over $Y$ into a field over $X$. One way to do this is to assemble, for every $x\in X$, all the fibres over points $y\in Y$ that satisfy $p(y) =x$ and to make a single fibre out of them. We call this  ``pushforward construction'' the \emph{forgetful map}.

In the first part of this section we introduce the forgetful map in a non-equivariant setting. The groupoids come back into play in the second part of the section. The non-equivariant construction can also be found in the book \cite{FellDo:88:1}, Paragraph 14.9; it is formulated in the language of Banach bundles rather than in the language of u.s.c.~fields of Banach spaces.

\subsubsection{The forgetful map for fields}

Let $X$ and $Y$ be locally compact Hausdorff spaces and let $p\colon Y\to X$ be continuous.

\begin{defprop}[The pushforward $p_*E$] Let $E$ be a u.s.c.~field of Banach spaces over $Y$. For all $x\in X$, define\footnote{This definition makes sense if $x\in p(Y)$, and can and should be interpreted as $p_*(E)_x=0$ if $x\notin p(Y)$.}
\[
p_*(E)_x:=\ContSect_0\left(Y_x,\ E\restr_{Y_x}\right).
\]
On this family of Banach spaces over $X$, define a structure of a u.s.c.~field of Banach spaces over $X$ as follows: For all $\xi\in \ContSect_0(Y,E)$, define $p_*(\xi)\colon x\mapsto \xi\restr_{Y_x}$. Then $\ContSect_0:= \left\{p_*(\xi):\ \xi \in \ContSect_0(Y,E)\right\}$ satisfies conditions (C1) - (C3) of the definition of a u.s.c.~field of Banach spaces and therefore defines\footnote{Use Proposition 1.1.4 of \cite{Lafforgue:06}.} a structure of a u.s.c.~field of Banach spaces over $X$ on $p_*(E)$. It has the property $\ContSect_0=\ContSect_0\left(X,\ p_*E\right)$.
\end{defprop}
\begin{proof}
Let $x\in X$. Then $\ContSect_0(Y_x, E\restr_{Y_x})$ is a $\C$-linear subspace of the space $\prod_{x\in X} p_*(E)_x$, i.e., it satisfies (C1). The space $Y_x$ is a closed subspace of $Y$, so we can apply Proposition~E.5.2 of \cite{Paravicini:07} to it which says that the map $\xi \mapsto \xi\restr_{Y_x}$ is a metric surjection from $\ContSect_0(Y,E)$ onto $\ContSect_0(Y_x,\ E\restr_{Y_x})$. In particular, the set $\ContSect_0$ satisfies (C2). As for (C3), let $x_0\in X$, $\varepsilon>0$ and $\xi \in \ContSect_0(Y,E)$. Let $L:= \left\{y\in Y:\ \norm{\xi(y)} \geq \norm{p_*(\xi)(x_0)} + \varepsilon\right\}$. Then $L$ is a compact subset of $Y$ because $\xi$ is vanishing at infinity. Hence its image $K:= p(L)$ is a compact subset of $X$. This set $K$ does not contain $x_0$, so its complement $U:= X\setminus K$ is an open neighbourhood of $x_0$ such that for $u\in U$ we have
\[
\norm{p(\xi)(u)} = \sup_{y\in Y_x} \norm{\xi(y)} \leq \norm{p_*(\xi)(x_0)} +\varepsilon,
\]
where the supremum is assumed to be zero if taken over the empty set. Hence we have shown that $\abs{p_*(\xi)}$ is upper semi-continuous.

It remains to show that $\ContSect_0=\ContSect_0(X,\ p_*E)$. Let $\xi$ be in $\ContSect_0(Y,E)$ and $\varepsilon>0$. Find a compact subset $L$ of $Y$ such that $\norm{\xi(y)} \leq \varepsilon$ whenever $y\in Y\setminus L$. Let $K:= p(L)$. Then $K$ is a compact subset such that $\norm{p_*(\xi)(x)} = \sup_{y\in Y_x} \norm{\xi(y)} \leq \varepsilon$ for all $x\in X\setminus K$. So $p_*(\xi)$ vanishes at infinity. This shows that $\xi \mapsto p_*\xi$ is an (isometric) map from $\ContSect_0(Y,\ E)$ to $\ContSect_0(X,\ p_*E)$. The image is invariant under multiplication with elements of $\Cont_c(X)$ and dense in all fibres, so it is dense (see, for example, Proposition~3.1.27 of \cite{Paravicini:07}). Hence the image is all of $\ContSect_0(X,\ p_*E)$.
\end{proof}

\noindent The pushforward construction can be extended to linear and bilinear maps:
\begin{itemize}
\item Let $E$ and $F$ be u.s.c.~fields of Banach spaces over $Y$ and let $T$ be a bounded continuous field of linear maps from $E$ to $F$. For all $x\in X$, define
\[
p_*(T)_x\colon p_*(E)_x \to p_*(F)_y,\ \xi\mapsto \left[Y_x\ni y\mapsto T_y(\xi(y))\right].
\]
Then $p_*T$ is a continuous field of linear maps bounded by $\norm{T}$.
\item Let $E_1$, $E_2$ and $F$ be u.s.c.~fields of Banach spaces over $Y$ and let $\mu\colon E_1 \times_Y E_2 \to F$ be a bounded continuous field of bilinear maps. For all $x\in X$, define
\[
p_*(\mu)_x \colon p_*(E_1)_x \times p_*(E_2)_x \to p_*(F)_x,\ (\xi_1,\xi_2) \mapsto \left[Y_x\ni y\mapsto \mu_y(\xi_1(y),\ \xi_2(y))\right].
\]
Then $p_*\mu$ is a continuous field of bilinear maps bounded by $\norm{\mu}$. If $\mu$ is non-degenerate, then so is $p_* \mu$, and vice versa.
\end{itemize}

\noindent Because this definition respects the associativity of bilinear maps, $p_*A$ is a u.s.c.~field of Banach algebras over $X$ if $A$ is a u.s.c.~field of Banach algebras over $Y$, canonically. The situation is similar for Banach modules and pairs.

We now state a technical result concerning the forgetful map; the proof can be found in \cite{Paravicini:07}, Section~8.3.1. This result is crucial if one wants to check all the details in the upcoming section.

\begin{proposition}[Pushforward and Pullback]\label{Proposition:FieldsPushForwardAndPullBack} Let $Y'$ be another locally compact Hausdorff space and let $p'\colon Y' \to X$ be continuous. Let $q\colon Y' \times_X Y \to Y$ and $q' \colon Y' \times_X Y\to Y'$ be the canonical ``projections''. Let $E$ be a u.s.c.~field of Banach spaces over $X$. Then $p'^*(p_*E) \cong q'_*(q^*E)$, i.e., the two ways of going from the upper right to the lower left corner in the following diagram yield the same result:
\[
\xymatrix{
Y' \times_X Y \ar[r]^-{q} \ar[d]_{q'} & Y\ar[d]^-{p} \\
Y' \ar[r]^{p'} & X
}
\]
\end{proposition}

\subsubsection{The forgetful map for equivariant fields}

For proofs of most of the results of this subsection we refer to \cite{Paravicini:07}, Section~8.3.2.

In this subsection, let $\gG$ be a locally compact Hausdorff groupoid and let $Y$ be a locally compact Hausdorff left $\gG$-space with anchor map $\rho\colon Y\to X=\gG^{(0)}$.

\begin{definition} Let $E$ be a $\gG \ltimes Y$-Banach space with action $\alpha$. Let $\gamma\in \gG$ and $\xi_{s(\gamma)} \in (\rho_*E)_{s(\gamma)} = \ContSect_0\left(Y_{s(\gamma)},\ E\restr_{Y_{s(\gamma)}}\right)$. Define a section $\gamma \xi_{s(\gamma)}\in \ContSect_0\left(Y_{r(\gamma)},\ E\restr_{Y_{r(\gamma)}}\right) = (\rho_*E)_{r(\gamma)}$ by
\[
\left(\gamma \xi_{s(\gamma)}\right)\ (y) := (\gamma, y) \left(\xi_{s(\gamma)} \left( \gamma^{-1} y\right)\right)
\]
for all $y\in Y_{r(\gamma)}$. This defines an action of $\gG$ on $\rho_* E$.
\end{definition}
\noindent This construction respects  the equivariance of linear and bilinear maps. Hence the pushforward of a $\gG\ltimes Y$-Banach algebra is a $\gG$-Banach algebra etc. For clarity, we detail the construction for equivariant Banach pairs:

Let $B$ be a $\gG \ltimes Y$-Banach algebra and let $E=(E^<,E^>)$ be a $\gG \ltimes Y$-Banach $B$-pair with bracket $\langle\cdot,\cdot\rangle_E$. Then $\rho_*E= \left(\rho_* E^<,\ \rho_* E^>\right)$ together with $\rho_* \langle\cdot,\cdot\rangle_E$ is a $\gG$-Banach $\rho_* B$-pair. This construction respects equivariant concurrent homomorphisms of Banach pairs. If $F$ is another right $\gG\ltimes Y$-Banach $B$-pair and $T\in \Lin_B(E,F)$ is a $B$-linear bounded continuous field of operators, then $\rho_* T$ is in $\Lin_{\rho_*B}\left(\rho_* E,\ \rho_* F\right)$ with $\norm{\rho_* T} \leq \norm{T}$. Note that, if we express the actions on $E$ and $\rho_* E$ in terms of the ``unitary'' operators, we have under suitable identifications using Proposition~\ref{Proposition:FieldsPushForwardAndPullBack}:
\begin{equation}\label{Equation:PushforwardAndUnitaryAction}
V_{\rho_* E} = \pi_{1,*} V_E
\end{equation}
where $\pi_1 \colon \gG \ltimes Y \to \gG$ is the projection on the first component; see \cite{Paravicini:07:Induction:arxiv}, Section~3.2 for more information about the operator $V_E$. 

The following lemma is straightforward to check.

\begin{lemma}
Let $B$ be a $\gG \ltimes Y$-Banach algebra and let $E$ and $F$ be $\gG \ltimes Y$-Banach $B$-pairs. Let $\xi^<\in \ContSect_0(X,E^<)$ and $\eta^>\in \ContSect_0(X,F^>)$. Then
\[
\rho_* \left(\ketbra{\eta^>}{\xi^<}\right) = \ketbra{\rho_*(\eta^>)}{\rho_*(\xi^<)} \in \Komp_{\rho_*B}\left(\rho_* E,\ \rho_* F\right).
\]
It follows that $\rho_*\left(\Komp_B(E,F)\right) \subseteq \Komp_{\rho_* B} \left(\rho_*E,\ \rho_* F\right)$.
\end{lemma}

\begin{proposition}
Let $A$ and $B$ be $\gG \ltimes Y$-Banach algebras. Let $(E,T) \in \EbanW{\gG \ltimes Y}\left(A,\ B\right)$. Then $\left(\rho_* E,\ \rho_* T\right)$ is in $\EbanW{\gG}\left(\rho_* A,\ \rho_* B\right)$.
\end{proposition}
\begin{proof}
Surely, $E$ is a graded non-degenerate $\gG$-Banach $\rho_*A$-$\rho_*B$-pair and $\rho_*T$ is an odd continuous field of linear operators on $\rho_*E$. Now let $a\in \ContSect_0\left(X,\ \rho_*A\right)$. Then there is an $a'\in \ContSect_0\left(Y,\ A\right)$ such that $a= \rho_* a'$. Now
\[
\left[a,\ \rho_*T\right] = \left[\rho_* a',\ \rho_* T\right] = \rho_* \left[a', T\right] \in \Komp_{\rho_*B}\left(\rho_*E\right).
\]
Similarly,
\[
a\left(\rho_*T^2-1\right) = \rho_* \left(a'(T^2-1)\right) \in \Komp_{\rho_*B}\left(\rho_*E\right).
\]
Now let $\tilde{a} \in \ContSect_0\left(\gG,\ r_{\gG}^*\rho_* A\right)$. As above, define $\pi_1\colon \gG\ltimes Y\to\gG,\ (\gamma,y)\mapsto \gamma$. We identify $r_{\gG}^* \rho_* A$ and $\pi_{1,*} r_{\gG\ltimes Y}^* A$ (and do the same for $B$ and $E$) and regard $\tilde{a}$ as an element of $\ContSect_0\left(\gG,\ \pi_{1,*} r_{\gG\ltimes Y}^*  A\right)$. We can then find an element $\tilde{a}'\in \ContSect_0\left(\gG\ltimes Y,\ r_{\gG \ltimes Y}^* A\right)$ such that $\pi_{1,*} \tilde{a}'= \tilde{a}$. Using Equation~(\ref{Equation:PushforwardAndUnitaryAction}) and suitable identifications we can now conclude
\begin{eqnarray*}
&& \tilde{a}\left(V_{\rho_*E} \circ (s_{\gG}^* \rho_*T)\circ V_{\rho_*E}^{-1} - r_{\gG}^*\rho_*T\right)\\
& =& \pi_{1,*} \tilde{a}' \left((\pi_{1,*} V_{E})\circ (\pi_{1,*} s^*_{\gG\ltimes Y}(T))\circ (\pi_{1,*} V_{E}^{-1}) - \pi_{1,*} r^*_{\gG\ltimes Y}(T)\right)\\
& =& \pi_{1,*} \left(\tilde{a}' \left(V_E \circ  s_{\gG \ltimes Y}^*T \circ V_E^{-1}- r_{\gG\ltimes Y}^* T\right)\right) \in \Komp_{\pi_{1,*}r_{\gG\ltimes Y}^*B}\left(\pi_{1,*}r_{\gG\ltimes Y}^* E\right) = \Komp_{r_{\gG}^*\rho_* B}\left(r_{\gG}^*\rho_* E\right).
\end{eqnarray*}
Proposition~3.13 of \cite{Paravicini:07:Induction:arxiv} tells us that we can define $\KKbanW{\gG}$-cycles between $\gG$-Banach algebras with locally compact Hausdorff $\gG$ also using compact instead of locally compact operators, so we have shown that $(\rho_* E, \rho_* T)$ is a $\KKbanW{\gG}$-cycle. 
\end{proof}

\noindent Let $B$ be a $\gG\ltimes Y$-Banach algebra. Then $\rho_*(B[0,1])$ is canonically isomorphic to $(\rho_*B)[0,1]$. Moreover, if $E$ is a left $\gG\ltimes Y$-Banach $B$-module and is $F$ a right $\gG\ltimes Y$-Banach $B$-module, one of them being non-degenerate, then
\[
\rho_*(E) \otimes_{\rho_*(B)} \rho_*(F) \cong \rho_*\left(E \otimes_B F\right).
\]
The proof of this result can be found in \cite{Paravicini:07}, Lemma 8.3.17; it uses a technical result which is published separately in \cite{Paravicini:07:CNullX:erschienen}.

We can now conclude, using some standard arguments, that the forgetful map lifts to the $\KKban$-groups:

\begin{theorem}\label{Theorem:ForgetfulHomomorphism}
Let $A$ and $B$ be $\gG \ltimes Y$-Banach algebras. Then $\rho_*$ gives a homomorphism
\[
\rho_* \colon \KKbanW{\gG \ltimes Y} \left(A,B\right) \to \KKbanW{\gG} \left(\rho_*A,\ \rho_*B\right).
\]
\end{theorem}

\subsection{The forgetful map and the descent}

In Paragraph~\ref{Subsection:ForgetfulAndAssembly}, we are going to show that the forgetful map of Theorem~\ref{Theorem:ForgetfulHomomorphism} is compatible with the Bost assembly map. Before we analyse that situation, we have to first discuss how the forgetful map and the descent are related. So let $\gG$ be a locally compact Hausdorff groupoid equipped with a Haar system. Let $Y$ be a locally compact Hausdorff $\gG$-space with anchor map $\rho$.

\begin{definition}\label{Definition:TheInjectionTildeIota} Let $E$ be a $\gG\ltimes Y$-Banach space. For all $\xi \in \ContSect_c\left(\gG \ltimes Y,\ r_{\gG \ltimes Y}^*  E\right)$, define $\hat{\iota}_E (\xi) \in \ContSect_c\left(\gG,\ r_{\gG}^*\rho_* E\right)$ by
\[
\hat{\iota}_{E} (\xi)\ (\gamma):= \left[Y_{r_{\gG}(\gamma)} \ni y \mapsto \xi(\gamma, y)\right] \in \rho_* E _{r_{\gG}(\gamma)}
\]
for all $\gamma \in \gG$. The map $\hat{\iota}_{E}$ is continuous for the inductive limit topologies, injective, and has dense image.
\end{definition}

Let $E_1$, $E_2$ and $F$ be $\gG\ltimes Y$-Banach spaces and let $\mu\colon E_1 \times_Y E_2  \to F$ be a bounded equivariant continuous field of bilinear maps. Then $(\rho_*\mu) \left(\hat{\iota}_{E_1}(\xi_1),\ \hat{\iota}_{E_2}(\xi_2)\right) = \hat{\iota}_{F}\left(\mu(\xi_1,\xi_2)\right)$ for all $\xi_1\in \ContSect_c\left(\gG\ltimes Y,\ r_{\gG \ltimes Y}^*E_1\right)$ and $\xi_2\in \ContSect_c\left(\gG\ltimes Y,\ r_{\gG \ltimes Y}^*E_2\right)$; this could also be written as
\[
\hat{\iota}_{E_1}(\xi_1)\ *\ \hat{\iota}_{E_2}(\xi_2) = \hat{\iota}_{F}\left(\xi_1 *\xi_2\right).
\]
In particular, this applies to the multiplication of a Banach algebra; more precisely, for every $\gG\ltimes Y$-Banach algebra $B$, the map $\hat{\iota}_{B}$ is a continuous injective homomorphism with dense image.

\begin{defprop}\label{DefProp:TildeIotaAndUnconditionalCompletions} Let $\mA(\gG)$ be an unconditional completion of $\Cont_c(\gG)$ and let $E$ be a $\gG \ltimes Y$-Banach space. For all $\xi \in \ContSect_c\left(\gG\ltimes Y,\ r_{\gG\ltimes Y}^* E\right)$, define
\[
\norm{\xi}_{\mA_Y}:= \norm{\hat{\iota}_{E} (\xi)}_{\mA}.
\]
This defines a semi-norm on $\ContSect_c\left(\gG\ltimes Y,\ r_{\gG\ltimes Y}^* E\right)$. If $E=\C_Y$, then $\ContSect_c\left(\gG\ltimes Y,\ r_{\gG\ltimes Y}^* \C_Y\right) = \Cont_c\left(\gG\ltimes Y\right)$ and $\norm{\cdot}_{\mA_Y}$ is an unconditional semi-norm on $\Cont_c\left(\gG\ltimes Y\right)$. The map $\hat{\iota}_E$ extends to an isomorphism on the completions:
\[
\hat{\iota}_{E}\colon \mA_Y\left(\gG\ltimes Y,\ E\right) \cong \mA\left(\gG,\ \rho_* E\right).
\]
If $B$ is a $\gG\ltimes Y$-Banach algebra, then
\[
\hat{\iota}_{B}\colon \mA_Y\left(\gG\ltimes Y,\ B\right) \cong \mA\left(\gG,\ \rho_* B\right)
\]
as Banach algebras.
\end{defprop}

\begin{proposition}\label{Proposition:PushforwardDescentCommute}
Let $A$ and $B$ be $\gG\ltimes Y$-Banach algebras and let $\mA(\gG)$ be an unconditional completion of $\gG$. Then the following diagram is commutative:
\[
\xymatrix{
\KKbanW{\gG\ltimes Y}\left(A,\ B\right) \ar[d]_{\rho_*} \ar[rr]^-{j_{\mA_Y}} && \KKban\left(\mA_Y\left(\gG\ltimes Y,A\right),\ \mA_Y\left(\gG\ltimes Y, B\right)\right)\ar[d]^{\cong} \\
\KKbanW{\gG}\left(\rho_*A,\ \rho_*B\right) \ar[rr]^-{j_{\mA}} && \KKban\left(\mA\left(\gG,\ \rho_*A\right),\ \mA\left(\gG,\ \rho_*B\right)\right)
}
\]
where the isomorphism on the right-hand side is given by the isomorphism of Proposition~\ref{DefProp:TildeIotaAndUnconditionalCompletions} in both variables.
\end{proposition}
\begin{proof}
This is true already on the level of cycles: Let $(E,T) \in \EbanW{\gG\ltimes Y} \left(A,\ B\right)$. Then $\rho_*(E,T) = \left(\rho_*E,\ \rho_*T\right)$ by definition. The module $\mA\left(\gG, \rho_*E^>\right)$ is a completion of $\ContSect_c\left(\gG,\ r_{\gG}^*\rho_* E^>\right)$ for the norm $\norm{\cdot}_{\mA}$. Using the continuous injective linear map $\hat{\iota}_{E^>}$ from $\ContSect_c\left(\gG\ltimes Y,\ r_{\gG\ltimes Y}^*E^>\right)$ to $\ContSect_c\left(\gG,\ r_{\gG}^*\rho_* E^>\right)$ introduced in \ref{Definition:TheInjectionTildeIota} we get a linear isometric isomorphism $\hat{\iota}_{E^>}\colon\mA_Y\left(\gG\ltimes Y,\ E^>\right) \to \mA\left(\gG,\ \rho_*E^>\right)$; analogously, we get a linear isomorphism $\hat{\iota}_{E^<}\colon\mA_Y\left(\gG\ltimes Y,\ E^<\right) \to \mA\left(\gG,\ \rho_*E^<\right)$. Together, this gives an isomorphism of Banach pairs $\hat{\iota}_E$ from $\mA_Y\left(\gG\ltimes Y,\ E\right)$ to $\mA\left(\gG, \rho_* E\right)$ with coefficient maps $\hat{\iota}_A$ and $\hat{\iota}_B$. It is straightforward to show that this isomorphism is compatible with the grading and intertwines the operators, i.e., it is an isomorphism of cycles.
\end{proof}

\subsection{The forgetful map and the assembly map}\label{Subsection:ForgetfulAndAssembly}

Let $\gG$ again be a locally compact Hausdorff groupoid with Haar system, and let $Y$ be a locally compact Hausdorff $\gG$-space with anchor map $\rho$. Let $B$ be a $\gG\ltimes Y$-Banach algebra. Let $Z$ be a locally compact $\gG \ltimes Y$-space; note that we can regard $Z$ also as a $\gG$-space. The space $Z$ is proper as a $\gG\ltimes Y$-space if and only if it is proper as a $\gG$-space. Moreover, $Z$ is $\gG \ltimes Y$-compact if and only if it is $\gG$-compact. We therefore get homomorphisms
\[
\KKbanW{\gG \ltimes Y}(\Cont_0(Z), B) \to \KKbanW{\gG} (\Cont_0(Z), \rho_*B) \to \KTh^{\top2,\ban}_0(\gG, \rho_*B),
\]
and, passing to the direct limit over $\gG\ltimes Y$-compact proper $Z$ and considering suspensions of $B$,
\[
\rho_* \colon \KTh^{\top2,\ban}_*(\gG \ltimes Y, B) \to \KTh^{\top2,\ban}_*(\gG, \rho_*B).
\]

Let $\mA(\gG)$ be an unconditional completion of $\Cont_c(\gG)$. Define the unconditional completion $\mA_Y(\gG \ltimes Y)$ of $\Cont_c(\gG\ltimes Y)$ as in \ref{DefProp:TildeIotaAndUnconditionalCompletions}. Recall that $\mA_Y(\gG \ltimes Y, B) \cong \mA(\gG ,\rho_*B)$.

\begin{proposition}\label{Proposition:ForgetfulHomomorphismAndAssemblyMap}
The following diagram is commutative
\begin{equation}\label{Diagram:ForgetfulMapAndAssemblyMap}
\xymatrix{
\KTh^{\top2,\ban}_*(\gG \ltimes Y, B) \ar[d] \ar[r] & \KTh_*(\mA_Y(\gG \ltimes Y, B))\ar[d] ^{\cong}\\
\KTh^{\top2,\ban}_*(\gG , \rho_* B) \ar[r] & \KTh_*(\mA(\gG, \rho_* B)).
}
\end{equation}
\end{proposition}

\begin{proof}
We show the following: If $Z$ is a locally compact $\gG \ltimes Y$-space, then the following diagram is commutative:
\[
\xymatrix{
\KKbanW{\gG\ltimes Y}(\Cont_0(Z),\ B) \ar[r] \ar[d] & \KKban(\mA_Y(\gG \ltimes Y, \Cont_0(Z)),\ \mA_Y(\gG\ltimes Y, B)) \ar[d]^{\cong} \ar[r]& \KTh_0(\mA_Y(\gG \ltimes Y, B))\ar[d]^{\cong}\\
\KKbanW{\gG}(\Cont_0(Z),\ \rho_* B) \ar[r]  & \KKban(\mA(\gG , \Cont_0(Z)),\ \mA(\gG, \rho_* B)) \ar[r]& \KTh_0(\mA(\gG , \rho_*B))
}
\]
The commutativity of the left square is a special case of Proposition~\ref{Proposition:PushforwardDescentCommute}.

To see that the right square commutes, one has to check the rather trivial fact that the two elements of $\KTh_0(\mA_Y(\gG \ltimes Y, \Cont_0(Z)))$ and $\KTh_0(\mA(\gG, \Cont_0(Z)))$ which are used to define the horizontal maps can be identified. Let $c$ be a cut-off function on $Z$ for $\gG$. Then $c$ is also a cut-off function for $\gG\ltimes Y$.

Define $p(\gamma, y, z):= c^{1/2} (z) c^{1/2}((\gamma^{-1},y)z)$ for all $(\gamma,y,z)\in (\gG\ltimes Y) \ltimes Z$.

We can regard $p$ as an element of $\ContSect_c(\gG\ltimes Y, r^* \Cont_0(Z))$, but also as an element of $\Cont_c(\gG, \Cont_0(Z))$ because $(\gG\ltimes Y) \ltimes Z$ is homeomorphic to $\gG\times Z$ (just forget the (redundant) entry in $Y$). The element of $\mA_{Y}(\gG\ltimes Y, \Cont_0(Z))$ given by $p$ corresponds to the element of $\mA(\gG,\Cont_0(Z))$ given by $p$. Hence the square commutes because the vertical arrows only stem from the identification of these two algebras.
\end{proof}

\begin{corollary}\label{Corollary:ForgetfulMapAndAssemblyMap:CStar}
Let $G$ be a locally compact Hausdorff second countable \emph{group} and $Y$ be a locally compact Hausdorff second countable $G$-space. Let $B$ be a separable $G\ltimes Y$-C$^*$-algebra. Then the following diagram is commutative
\begin{equation}\label{Diagram:ForgetfulMapAndAssemblyMap:CStar}
\xymatrix{
\KTh^{\top2}_*(G \ltimes Y, B) \ar[d]_{\cong} \ar[r] & \KTh_*(\mA_Y(G \ltimes Y, B))\ar[d] ^{\cong}\\
\KTh^{\top2}_*(G , \ContSect_0(Y, B)) \ar[r] & \KTh_*(\mA(G, \ContSect_0(Y, B))).
}
\end{equation}
In particular, the Bost conjecture with separable C$^*$-algebra coefficients is true for $G\ltimes Y$, $\mA_Y$ and $B$ if and only if it is true for $G$, $\mA$ and $\ContSect_0(Y,B)$.
\end{corollary}
\begin{proof}
First note that $\rho_*B$ is $\ContSect_0(Y,B)$, because the unit space of $G$ is a point.

What we have to check is that we can replace $\KTh^{\top2,\ban}_*$ with $\KTh^{\top2}_*$ in the diagram. But this is straightforward because it is easy to show that the canonical map from $\KTh^{\top2}_*$ to $\KTh^{\top2,\ban}_*$ is compatible with the forgetful map.

As mentioned above, it is proved in \cite{CEO:03} that the forgetful map for C$^*$-algebras is an isomorphism. Hence the corollary follows.
\end{proof}

\section{The special case of $\Leb^1$-algebras} \label{Section:TheSpecialCaseOfL1}

We now analyse the right-hand side of Diagram~(\ref{Diagram:ForgetfulMapAndAssemblyMap}) of Proposition~\ref{Proposition:ForgetfulHomomorphismAndAssemblyMap} closer in the case that $\mA(\gG)$ is $\Leb^1(\gG)$. We confine ourselves to the case that the groupoid $\gG$ is actually a group $G$. Let $Y$ be a locally compact Hausdorff $G$-space (and $\rho$ be the constant map to the one-point unit space of the ``groupoid'' $G$).

The important technical detail that complicates matters is that there is a difference between the unconditional completions $\Leb^1(G \ltimes Y)$ and $\Leb^1_Y(G \ltimes Y)$ of $\Cont_c(G\ltimes Y)$. By $\Leb^1(G\ltimes Y)$ we mean the usual unsymmetrised $\Leb^1$-algebra of the groupoid $G\ltimes Y$, so its norm is given by
\[
\norm{\xi}_{\Leb^1}= \sup_{y \in Y} \int_{g \in G}  \abs{\xi(g,y)} \rmd g, \quad \xi \in \Cont_c(G\ltimes Y).
\]
The algebra $\Leb^1_Y(G \ltimes Y)$ is defined by \ref{DefProp:TildeIotaAndUnconditionalCompletions}, so it satisfies $\Leb^1_Y(G\ltimes Y) \cong \Leb^1(G, \Cont_0(Y))$, canonically, and hence the norm is given by
\[
\norm{\xi}_{\Leb^1_Y}=\int_{g \in G}\ \sup_{y \in Y} \abs{\xi(g,y)} \rmd g, \quad \xi \in \Cont_c(G\ltimes X).
\]
It follows that $\norm{\xi}_{\Leb^1}\leq \norm{\xi}_{\Leb^1_Y}$ for all $\xi\in \Cont_c(G \ltimes Y)$, but in general, equality does not hold.

If $B$ is a $G\ltimes Y$-Banach algebra, then we nevertheless get a norm-decreasing homomorphism
\[
\psi_B\colon \Leb^1(G, \ContSect_0(Y,B)) \cong \Leb^1_Y(G\ltimes Y, B) \to \Leb^1(G\ltimes Y, B);
\]
note that $\rho_*B = \ContSect_0(X,B)$. By Proposition~\ref{Proposition:ForgetfulHomomorphismAndAssemblyMap} and \cite{Paravicini:07:Induction:arxiv}, Proposition~6.3, we thus obtain a commutative diagram 
\[
\xymatrix{
\KTh^{\top2,\ban}_*(G \ltimes Y, B) \ar[d] \ar[r] & \KTh_*(\Leb^1(G \ltimes Y, B))\\
\KTh^{\top2,\ban}_*(G , \ContSect_0(Y,B)) \ar[r] & \KTh_*(\Leb^1(G, \ContSect_0(Y,B)))\ar[u]_{\psi_{B,*}}.
}
\]
Although the two algebras $\Leb^1(G, \ContSect_0(Y,B))$ and $\Leb^1(G\ltimes Y, B)$ might not be isomorphic, it still makes sense to ask whether $\psi_B$ is an isomorphism in $\KTh$-theory. We have the following positive result in this direction:

\begin{proposition}\label{Proposition:L1CompletionsIsomorphicKTheory}
If $Y$ is a proper or a discrete $G$-space, then $\psi_{B,*}$ is an isomorphism:
\[
\KTh_*(\Leb^1(G, \ContSect_0(Y,B))) \cong \KTh_*(\Leb^1(G \ltimes Y, B)).
\]
\end{proposition}
\begin{proof}
The method that we use in both cases is to produce a dense hereditary subalgebra $A$ of $\Leb^1_Y(G\ltimes Y, B) \cong \Leb^1(G, \ContSect_0(Y,B))$ which is also (dense and) hereditary in $\Leb^1(G\ltimes Y, B)$. See Lemme~1.7.8 and Lemme~1.7.10 of \cite{Lafforgue:02}.

In the case that $Y$ is proper, a good choice for such an algebra is $A:=\ContSect_c(G\ltimes Y, r^*B)$, where $r:G\ltimes Y \to Y$ is the range map of the groupoid $G\ltimes Y$. Recall that both, $\Leb^1(G \ltimes Y, B)$ and $\Leb^1_Y(G\ltimes Y, B) \cong \Leb^1(G, \ContSect_0(Y,B))$ are actually defined as completions of $A$.

If $\xi_1$ and $\xi_2$ are elements of $A = \ContSect_c(G\ltimes Y, r^*B)$. If $\xi$ is another element of $A$, then the support of $\xi_1 * \xi * \xi_2$ is contained in the set $K:= \{ \gamma \in G\ltimes Y:\ r(\gamma) \in r(\supp \xi_1) \wedge s(\gamma) \in (\supp\xi_2)\}$. This set is compact. Secondly, it is easy to see that
\[
\norm{\xi_1 * \xi * \xi_2}_{\infty} \leq \norm{\xi_1}_{\Leb^1} \norm{\xi}_{\Leb^1} \norm{\xi_2}_{\infty}.
\]
It follows that $\xi_1 * \Leb^1(G\ltimes Y, B) * \xi_2$ is contained in $A=\ContSect_c(G\ltimes Y, r^*B)$. A fortiori, $A$ is hereditary in $\Leb^1(G\ltimes Y, B)$. Compare the proof of Lemme~1.7.8 in \cite{Lafforgue:02}.

Note that this choice for $A$ seems to work only if $Y$ is proper.

However, if $Y$ is discrete, we can construct an algebra which is larger than $\ContSect_c(G\ltimes Y, r^* B)$ (and hence dense in both completions), but which is also hereditary. It is not clear whether a similar idea works for general $Y$. Assume from now on that $Y$ is discrete.

For each subset $M$ of $Y$, define
\[
A_M:= \{\xi \in \ContSect_c(G\ltimes Y, r^*B ):\ r(\supp \xi) \subseteq M\}.
\]
Note that $A_M$ is a right ideal of $\ContSect_c(G\ltimes Y, r^*B)$.

\begin{lemma*}
If $M$ is a finite set, then the norm on $A_M$ which is inherited from $\Leb^1(G\ltimes Y, B)$ is equivalent to the norm on $A_M$ which is inherited from $\Leb^1(G, \ContSect_0(Y,B))$.
\end{lemma*}
\begin{proof}
Let $\xi \in A_M$. Then
\begin{eqnarray*}
\norm{\xi}_{\Leb^1} &\leq& \norm{\xi}_{\Leb^1_Y} = \int_{g\in G} \sup_{y \in Y} \norm{\xi(g,y)} \rmd g \\
&=& \int_{g\in G} \sup_{m \in M} \norm{\xi(g,m)} \rmd g \leq \int_{g\in G} \sum_{m \in M} \norm{\xi(g,m)} \rmd g\\
&=&  \sum_{m \in M} \int_{g\in G} \norm{\xi(g,m)} \rmd g \leq \abs{M} \sup_{m \in M} \int_{g\in G} \norm{\xi(g,m)} \rmd g\\
&=& \abs{M} \norm{\xi}_{\Leb^1}.
\end{eqnarray*}
\end{proof}

\noindent It follows from this lemma that $\psi$ maps the closure $\overline{A}_M$ of $A_M$ in $\Leb^1(G, \ContSect_0(Y,B))$ bijectively and bicontinuously onto the closure of $A_M$ in $\Leb^1(G\ltimes Y, B)$ (which we also call $\overline{A}_M$). Note that $\overline{A}_M$ is a right ideal of both, $\Leb^1(G, \ContSect_0(Y,B))$ and $\Leb^1(G\ltimes Y, B)$, because $A_M$ is a right ideal of $\ContSect_c(G\ltimes Y, r^*B)$.

If $M$ and $N$ are finite subsets of $Y$ with $M\subseteq N$, then $A_M \subseteq A_N$ and hence $\overline{A}_M \subseteq \overline{A}_N$.

Define
\[
A:= \bigoplus_{M \subseteq Y \text{ finite}} \overline{A}_M \subseteq \Leb^1(G,\ContSect_0(Y,B)).
\]
Note that $\psi$ is injective on $A$, so we can think of $A$ as a subspace also of $\Leb^1(G\ltimes Y, B)$. Indeed, $A$ is a linear subspace of both algebras, and because all the $\overline{A}_M$ are right ideals, also $A$ is a right ideal (in both completions). Moreover, the union of all $A_M$ contains $\ContSect_c(G\ltimes Y, r^*B)$, so $A$ is dense in $\Leb^1(G, \ContSect_0(Y,B))$ and $\Leb^1(G\ltimes Y, B)$. In particular, $A$ is a dense hereditary subalgebra of both Banach algebras. By Lemme~1.7.10 in \cite{Lafforgue:02}, this means that $\psi$ is an isomorphism in $\KTh$-theory.
\end{proof}

\begin{remark}
If $Y$ is not discrete, then the preceding argument is no longer valid. The crucial point is that, for the argument to work, the $\Leb^1$-norm and the $\Leb^{\infty}$-norm for the continuous functions on any given compact subset of $Y$ have to be equivalent; this is a very restrictive condition on $Y$.
An instructive counter-example is $G=\R$ and $Y= S^1$ with the action by rotation. It is easy to construct a continuous function $f$ on $G\times Y = \R \times S^1$ which vanishes at infinity such that $\sup_{y\in S^1} \int_{t\in \R} \abs{f(t,y)} \rmd t < \infty$ and $\int_{t\in \R} \sup_{y\in S^1} \abs{f(t,y)} \rmd t = \infty$.

However, there might be a different way to construct a dense hereditary subalgebra of $\Leb^1(G\ltimes Y)$ which is also contained in $\Leb^1(G, \Cont_0(Y))$.
\end{remark}

If $Y$ is discrete, we can now improve Corollary~\ref{Corollary:ForgetfulMapAndAssemblyMap:CStar} for $\Leb^1$-algebras as follows:

\begin{corollary}\label{Corollary:ForgetfulMapAndL1AssemblyMap:CStar}
Let $G$ be a locally compact Hausdorff second countable \emph{group} and let $Y$ be a locally compact Hausdorff second countable $G$-space. Assume that $Y$ is discrete (or proper). Let $B$ be a separable $G\ltimes Y$-C$^*$-algebra. Then the following diagram is commutative
\begin{equation}\label{Diagram:ForgetfulMapAndL1AssemblyMap:CStar}
\xymatrix{
\KTh^{\top2}_*(G \ltimes Y, B) \ar[d]_{\cong} \ar[r] & \KTh_*(\Leb^1(G \ltimes Y, B))\\
\KTh^{\top2}_*(G , \ContSect_0(Y, B)) \ar[r] & \KTh_*(\Leb^1(G, \ContSect_0(Y, B)))\ar[u]_{\cong}.
}
\end{equation}
In particular, the $\Leb^1$-version of the Bost conjecture with separable C$^*$-algebra coefficients is true for $G\ltimes Y$ and $B$ if and only if it is true for $G$ and $\ContSect_0(Y,B)$.
\end{corollary}

\noindent Note that the proper case was already settled in \cite{Lafforgue:02} in the much stronger sense that the Bost conjecture is actually true for proper C$^*$-algebras.

\section{The Bost conjecture and open subgroups}\label{Section:Subgroups}

Let $G$ be a locally compact Hausdorff group and let $H$ be a closed subgroup of $G$. Let $B$ be an $H$-Banach algebra. Let $\Leb^1(H)$ and $\Leb^1(G)$ denote the (unsymmetrised) $\Leb^1$-completions of $\Cont_c(H)$ and $\Cont_c(G)$, respectively.

We define $\Ind_H^G B$ to be the following $G$-Banach algebra: The underlying Banach algebra is
\[
\Ind_H^G B:= \{f\colon G \to B:\ \forall g\in G\ \forall h\in H:\ f(gh^{-1}) = h f(g) \quad \text{and} \quad  \norm{f(g)} \to 0 \text{ if } gH \to \infty\};
\]
the $G$-action on $\Ind_H^G B$ is given by
\[
(gf)(g') = f(g^{-1}g')
\]
for all $g,g'\in G$ and $f\in \Ind_H^G B$. As for every $G$-Banach algebra, we can form the convolution algebra $\Leb^1(G, \Ind_H^G B)$.

The locally compact groupoid $G \ltimes G/H$ is equivalent to the group $H$ (as a groupoid) with equivalence $_{G\ltimes G/H}G_H$. Hence we get not only the standard induced $G$-Banach algebra $\Ind_H^G B$, but also the induced $G\ltimes G/H$-Banach algebra $\Ind_{H}^{G\ltimes G/H} B$, the main difference being that we keep track of the fibration over $G/H$ in the second case. If we write $\rho$ for the constant function from $G/H$ to the one-point set $\{e_G\}$ (which is the unit space of the groupoid $G$), then the forgetful functor $\rho_*$ from $G\ltimes G/H$-Banach algebras to $G$-Banach algebras takes $\Ind_{H}^{G\ltimes G/H}B$ to $\Ind_H^G B$:
\[
\rho_*(\Ind_H^{G\ltimes G/H} B) = \ContSect_0(G/H, \Ind_H^{G\ltimes G/H} B) = \Ind_H^G B.
\]
Assume now that $G/H$ is $\sigma$-compact; this is necessary because we want to apply the induction theorem of \cite{Paravicini:07:Induction:arxiv} which requires that the involved groupoids have $\sigma$-compact unit spaces.

We have the following commutative diagram:
\[
\xymatrix{
\KTh^{\top2,\ban}_*(H , B) \ar[d]_{\cong}\ar[rr]&&\KTh_*(\Leb^1(H,B))\ar[d]^{\cong} \\
\KTh^{\top2,\ban}_*(G\ltimes G/H, \Ind_{H}^{G\ltimes G/H} B)\ar[d]_{\rho_*}\ar[rr]&& \KTh_*(\Leb^1(G\ltimes G/H, \Ind_H^{G\ltimes G/H} B)) \\
\KTh^{\top2,\ban}_*(G,  \Ind_H^G B)\ar[rr] &&\KTh_*(\Leb^1(G, \Ind_H^{G} B))\ar[u]^{\psi_*}
}
\]
The upper square is just a special case of the induction theorem (\cite{Paravicini:07:Induction:arxiv}, Theorem~6.5). The lower square has been discussed in Section~\ref{Section:TheSpecialCaseOfL1}; here we write $\psi$ for the canonical homomorphism from $\Leb^1(G, \Ind_H^{G} B)$ to $\Leb^1(G\ltimes G/H, \Ind_H^{G\ltimes G/H} B)$. 

If $H$ is open (and $G/H$ is hence discrete) or $H$ is compact (so $G/H$ is proper), then the arrow in the above diagram labelled $\psi_*$ is an isomorphism by Proposition~\ref{Proposition:L1CompletionsIsomorphicKTheory}.

We can read the following corollaries off this diagram:

\begin{corollary}
The $\Leb^1$-version of the Bost conjecture with Banach algebra coefficients is true for $H$ and $B$ if and only if it is true for $G\ltimes G/H$ and $\Ind_{H}^{G\ltimes G/H} B$.
\end{corollary}

\begin{corollary}
The $\Leb^1$-version of the Bost conjecture with arbitrary Banach algebra coefficients is true for $H$ if and only if it is true for $G\ltimes G/H$.
\end{corollary}

\begin{corollary}
If $H$ is an open or compact subgroup of $G$, then $\KTh_*(\Leb^1(G, \Ind_H^G B))$ is isomorphic to $\KTh_*(\Leb^1(H,B))$.
\end{corollary}

The first two of these corollaries can be read off the upper square of the diagram, the third corollary follows from the right-hand part. As long as we don't have more information on the homomorphism $\rho_*$, we cannot say much more about the relation between the Bost conjecture with Banach algebra coefficients for $H$ and for $G$.

The situation is different if $B$ is a $G$-C$^*$-algebra. In this case we can replace $\KTh^{\top2, \ban}_*$ with $\KTh^{\top2}_*$ in the above diagram. This is possible, because the canonical homomorphisms from $\KTh^{\top2}_*$ to $\KTh^{\top2,\ban}_*$ are compatible with induction and the forgetful map and because the C$^*$-version of the Bost assembly map factors through $\KTh^{\top2, \ban}$. We hence get the following commutative diagram:
\[
\xymatrix{
\KTh^{\top2}_*(H , B) \ar[d]_{\cong}\ar[rr]&&\KTh_*(\Leb^1(H,B))\ar[d]^{\cong} \\
\KTh^{\top2}_*(G\ltimes G/H, \Ind_{H}^{G\ltimes G/H} B)\ar[d]_{\rho_*}\ar[rr]&& \KTh_*(\Leb^1(G\ltimes G/H, \Ind_H^{G\ltimes G/H} B)) \\
\KTh^{\top2}_*(G,  \Ind_H^G B)\ar[rr] &&\KTh_*(\Leb^1(G, \Ind_H^{G} B))\ar[u]^{\psi_*}
}
\]
Recall that the forgetful map $\rho_*$ is an isomorphism in this case if $B$ is separable and $G$ is second countable. This is the main result of \cite{CEO:03}. We can hence deduce the following corollaries:

\begin{corollary}
The $\Leb^1$-version of the Bost conjecture with C$^*$-algebra coefficients is true for $H$ and $B$ if and only if it is true for $G\ltimes G/H$ and $\Ind_{H}^{G\ltimes G/H} B$.
\end{corollary}

\begin{corollary}
The $\Leb^1$-version of the Bost conjecture with arbitrary C$^*$-algebra coefficients is true for $H$ if and only if it is true for $G\ltimes G/H$.
\end{corollary}

\begin{corollary}
If $H$ is an open or compact subgroup of $G$ and $G$ is second countable, then not only is $\KTh_*(\Leb^1(G, \Ind_H^G B))$ isomorphic to $\KTh_*(\Leb^1(H,B))$, but also the following holds:
The $\Leb^1$-version of the Bost conjecture with C$^*$-algebra coefficients is true for $H$ and separable $B$ if and only if it is true for $G$ and $\Ind_H^G B$.
\end{corollary}

\noindent Finally, this directly implies the following theorem (which we only state for open subgroups as the case of compact subgroups is not interesting):

\begin{theorem}\label{Theorem:BostForOpenSubgroups}
Let $G$ be a second countable locally compact Hausdorff group. If the $\Leb^1$-version of the Bost conjecture with arbitrary separable C$^*$-algebra coefficients is true for $G$, then it is also true for every open subgroup of $G$.
\end{theorem}

\noindent In particular, this applies to countable discrete groups:

\begin{corollary}
Let $G$ a countable (discrete) group. If the $\Leb^1$-version of the Bost conjecture with arbitrary separable C$^*$-algebra coefficients is true for $G$, then it is also true for every subgroup of $G$.
\end{corollary}

\section{Group actions on trees}\label{Section:Tree}

Before we come to group actions on trees, we first analyse group actions on general discrete spaces.

\subsection{Actions on a discrete space}

Let $G$ be a locally compact Hausdorff group acting continuously on a discrete space $Z$. For all $z\in Z$ let $G_z$ and $Z_z$ be the stabiliser and the orbit of $z$, respectively. Let $S$ be a fundamental domain for the $G$-action on $Z$, so $Z = \coprod_{s\in S} Z_s$.

For each $s \in S$, the group $G_s$ and the groupoid $G \ltimes G/G_s$ are equivalent as groupoids, the equivalence being given by $G$. Note that the $G_s$ are open subgroups, so that we can apply the results of the preceding paragraph.

%

If $B$ is a $G$-Banach algebra and $z\in Z$, then we write $B_z$ for the Banach algebra $B$, equipped with the restriction of the $G$-action to $G_z$. We now show that we have an isomorphism of $G$-Banach algebras
\[
\Ind_{G_s}^{G} B_s\ \cong \ \Cont_0(Z_s, B)
\]
via the following map: If $f\in \Cont_0(Z_s, B)$, then define
\[
(\Phi(f))(g):= g^{-1} (f(g s))
\]
for all $g\in G$. It is straightforward to check that this defines a $G$-equivariant homomorphism $\Phi\colon \Cont_0(Z_s, B) \to \Ind_{G_s}^{G} B_s$. It is clearly injective, even isometric, and if $f\in \Ind_{G_s}^{G} B_s$, then we can define a map $f'$ from $Z_s$ to $B$ by choosing, for every $z\in Z_s$, a $g_z\in G$ such that $z= g_z s$; we can then define $f'(z) := g_z f(g_z)$. It is easy to see that this definition is actually independent of the choice of $g_z$ and $\Phi(f')=f$, so $\Phi$ is surjective.

%

Assume now that $B$ is not only a Banach algebra but a separable C$^*$-algebra. Then $\Phi$ is actually a $*$-isomorphism. Assume moreover that $G$ is second countable (from this it follows that the discrete spaces $G/G_s$ and $Z_s$ are (second) countable for all $s\in S$). Then we can use $\Phi$ and the results of the preceding paragraph to obtain, for every $s\in S$, a commutative diagram of the following form.
\[
\xymatrix{
\KTh^{\top2}_*(G_s, B_s) \ar[d]_{\cong} \ar[r] & \KTh_*(\Leb^1(G_s, B_s))\ar[d]^{\cong}\\
\KTh^{\top2}_*(G , \Cont_0(Z_s, B)) \ar[r] & \KTh_*(\Leb^1(G, \Cont_0(Z_s, B))).
}
\]
Taking the direct sum over all $s\in S$ we thus arrive at the following diagram
\[
\xymatrix{
\bigoplus_{s\in S} \KTh^{\top2}_*(G_s, B_s) \ar[d]_{\cong} \ar[r] & \bigoplus_{s\in S} \KTh_*(\Leb^1(G_s, B_s))\ar[d]^{\cong}\\
\bigoplus_{s\in S}\KTh^{\top2}_*(G , \Cont_0(Z_s, B)) \ar[r] \ar[d]_{\cong} & \bigoplus_{s\in S} \KTh_*(\Leb^1(G, \Cont_0(Z_s, B)))\ar[d]^{\cong}\\
\KTh^{\top2}_*(G , \Cont_0(Z, B)) \ar[r] & \KTh_*(\Leb^1(G, \Cont_0(Z, B))).
}
\]
Here we use the fact that $\bigoplus_{s\in S} \Cont_0(Z_s, B) \cong \Cont_0(Z,B)$ as C$^*$-algebras. That the vertical map in the lower left-hand corner is an isomorphism follows from the non-trivial fact that $\KTh^{\top2}$ is continuous for direct limits (see \cite{ChaEch:01:Permanence}, Section 7). That the vertical map in the lower right-hand corner is an isomorphism follows from the continuity of $\KTh$-theory for Banach algebras.

We can read the following result off the above diagram:

\begin{lemma}\label{Lemma:Stabilisers}
Let the second countable locally compact Hausdorff group $G$ act on a discrete space $Z$.
\begin{enumerate}
\item Let $B$ be a separable $G$-C$^*$-algebra. Then $G$ satisfies the $\Leb^1$-version of the Bost conjecture with coefficients in $\Cont_0(Z,B)$ if and only if the stabiliser $G_s$ satisfies the $\Leb^1$-version of the Bost conjecture with coefficients in $B_s$ for every $s\in S$.
\item Let all stabilisers satisfy the $\Leb^1$-version of the Bost conjecture with arbitrary separable C$^*$-algebra coefficients. Then $G$ satisfies the conjecture for the $G$-C$^*$-algebra $\Cont_0(Z,B)$ for every separable $G$-C$^*$-algebra $B$.
\end{enumerate}
\end{lemma}

\subsection{Actions on a tree}

\begin{theorem}\label{Theorem:Tree:MainTheorem}
Let $G$ be a locally compact second countable Hausdorff group which acts on an oriented tree $X$. Then the group $G$ satisfies the $\Leb^1$-version of the Bost conjecture with separable C$^*$-algebra coefficients if and only if the stabilisers of all vertices of $X$ satisfy it.
\end{theorem}
\begin{proof}
The only-if-part is a consequence of Theorem~\ref{Theorem:BostForOpenSubgroups} because the stabilisers are open subgroups of $G$.

To prove the if-part, we adjust the method detailed in \cite{OyonoOyono:97} for the Baum-Connes conjecture to the Bost conjecture; we remark that it follows directly from the results of \cite{ChaEch:01:Permanence} that the method can be generalised from countable discrete groups to arbitrary second countable groups. We hence introduce two six-term exact sequences (\ref{Equation:Tree:SixTermUp}) and (\ref{Equation:Tree:SixTermDown}) which are connected by Bost assembly maps.

Let $X^0$ and $X^1$ denote the spaces of vertices and edges of $X$, respectively. For every vertex $q\in X^0$ we write $X^1_q$ for the set of edges of $X$ that point towards $q$. Let
\[
\overline{X^1}:= X^1 \cup \{-\infty,+\infty\}
\]
be the following compactification of $X^1$: the discrete space $X^1$ is contained as an open subspace, the neighbourhoods of $+\infty$ are finite sections of sets of the form $\{+\infty\} \cup X^1_q$ and the neighbourhoods of $-\infty$ are finite sections of the complements of sets of the form $\{+\infty\} \cup X^1_q$; see also \cite{Pimsner:86}, \S 1. One can visualise $+\infty$ as a ``universal source'' of the tree and $-\infty$ as a ``universal target''. The situation is symmetric in the sense that, if we reverse the orientation of the tree $X$, then we get the same space with the r\^{o}les of $+\infty$ and $-\infty$ interchanged.

If $B$ is a $G$-Banach algebra, then we write $\Cont^+(X^1,B)$ for $\{f\in \Cont(\overline{X^1}, B):\ f(-\infty)=0\}$; this algebra fits into the exact sequence of $G$-Banach algebras
\begin{equation}\label{Equation:TreeExactSequence}
\xymatrix{ 0 \ar[r] & \Cont_0(X^1, B) \ar[r] & \Cont^+(X^1, B) \ar[r] & B \ar[r] & 0.
}
\end{equation}
The C$^*$-algebras $\Cont_0(X^0)$ and $\Cont^+(X^1)$ are equivariantly $\KK$-equivalent; this is Th\'{e}or\`{e}me~6.3.1 of \cite{OyonoOyono:97} which is a refinement of Proposition~14 of \cite{Pimsner:86}. If $B$ is a $G$-C$^*$-algebra, then it follows that $\Cont_0(X^0, B)$ and $\Cont^+(X^1, B)$ are equivariantly $\KK$-equivalent. It follows from \cite{ChaEch:01:Permanence}, Section 4, that there is a six-term exact sequence
\[
\xymatrix{
\KTh^{\top2}_0(G, \Cont_0(X^1, B)) \ar[r] & \KTh^{\top2}_0(G, \Cont^+(X^1, B)) \ar[r] & \KTh^{\top2}_0(G, B) \ar[d]\\
\KTh^{\top2}_1(G, B) \ar[u] & \KTh^{\top2}_1(G, \Cont^+(X^1, B)) \ar[l] & \KTh^{\top2}_1(G, \Cont_0(X^1, B)) \ar[l]
}
\]
in which we can replace $\Cont^+(X^1, B)$ with $\Cont_0(X^0, B)$ using the $\KK$-equivalence to obtain an exact sequence
\begin{equation}\label{Equation:Tree:SixTermUp}
\xymatrix{
\KTh^{\top2}_0(G, \Cont_0(X^1, B)) \ar[r] & \KTh^{\top2}_0(G, \Cont_0(X^0, B)) \ar[r] & \KTh^{\top2}_0(G, B) \ar[d]\\
\KTh^{\top2}_1(G, B) \ar[u] & \KTh^{\top2}_1(G, \Cont_0(X^0, B)) \ar[l] & \KTh^{\top2}_1(G, \Cont_0(X^1, B)) \ar[l]
}
\end{equation}

\noindent For discrete groups, this sequence appeared in \cite{OyonoOyono:97} where it is related to the corresponding sequence for the reduced crossed product algebras using the respective Baum-Connes assembly maps. We translate this to the setting of the Bost conjecture. To this end, let $\mA(G)$ be an unconditional completion of $\Cont_c(G)$. From (\ref{Equation:TreeExactSequence}) we obtain an exact sequence of Banach algebras
\[
\xymatrix{ 0 \ar[r] & \mA(G, \Cont_0(X^1, B)) \ar[r] & \mA(G, \Cont^+(X^1, B)) \ar[r] & \mA(G,B) \ar[r] & 0.
}
\]
It induces an exact sequence
\[
\xymatrix{
\KTh_0(\mA(G, \Cont_0(X^1, B))) \ar[r] & \KTh_0(\mA(G, \Cont^+(X^1, B))) \ar[r] & \KTh_0(\mA(G, B)) \ar[d]\\
\KTh_1(\mA(G, B)) \ar[u] & \KTh_1(\mA(G, \Cont^+(X^1, B))) \ar[l] & \KTh_1(\mA(G, \Cont_0(X^1, B))) \ar[l]
}
\]
Note that the preceding two exact sequences also work for general $G$-Banach algebras, but we prefer to work with $G$-C$^*$-algebras because we neither have an analogue of (\ref{Equation:Tree:SixTermUp}) for Banach algebras yet, nor can we generalise the following arguments to Banach algebras: Because equivariant $\KK$-equivalences of C$^*$-algebras descent to isomorphisms in $\KTh$-theory also in the setting of unconditional completions (Proposition~1.6.10 of \cite{Lafforgue:02}), we obtain an exact sequence
\begin{equation}\label{Equation:Tree:SixTermDown}
\xymatrix{
\KTh_0(\mA(G, \Cont_0(X^1, B))) \ar[r] & \KTh_0(\mA(G, \Cont_0(X^0, B))) \ar[r] & \KTh_0(\mA(G, B)) \ar[d]\\
\KTh_1(\mA(G, B)) \ar[u] & \KTh_1(\mA(G, \Cont_0(X^0, B))) \ar[l] & \KTh_1(\mA(G, \Cont_0(X^1, B))) \ar[l]
}
\end{equation}
As in \cite{OyonoOyono:97}, the Bost assembly maps connect the six-term exact sequences (\ref{Equation:Tree:SixTermUp}) and (\ref{Equation:Tree:SixTermDown}); this can be proved by adjusting the proof of Proposition~4.1 of \cite{ChaEch:01:Permanence} to our situation.

Assume now that $\mA(G) = \Leb^1(G)$ and that all the stabilisers of vertices of $X$ satisfy the $\Leb^1$-version of the Bost conjecture with separable C$^*$-algebra coefficients. Because the stabilisers of edges of $X$ are open subgroups of the stabilisers of vertices of $X$, the stabilisers of edges also satisfy the $\Leb^1$-version of the Bost conjecture with separable C$^*$-algebra coefficients.

Assume that $B$ is a separable $G$-C$^*$-algebra. If we consider the action of $G$ on $X^0$ and $X^1$, then we can deduce from Lemma~\ref{Lemma:Stabilisers} that $G$ satisfies the Bost conjecture with coefficient algebras $\Cont_0(X^0,B)$ and $\Cont_0(X^1, B)$. By the five-lemma, applied to the diagrams (\ref{Equation:Tree:SixTermUp}) and (\ref{Equation:Tree:SixTermDown}) and the Bost assembly maps connecting them, $G$ satisfies the $\Leb^1$-version of the Bost conjecture also with coefficients $B$.
\end{proof}

\nocite{ChaEchMey:01} \nocite{Dixmier:64} \nocite{DuGil:83} \nocite{EKQR:06} \nocite{JenThom:91} \nocite{MuhReWill:87} \nocite{Renault:80}


\begin{thebibliography}{EKQR06}

\bibitem[CE01]{ChaEch:01:Permanence}
J{\'e}r{\^o}me Chabert and Siegfried Echterhoff.
\newblock Permanence properties of the {B}aum-{C}onnes conjecture.
\newblock {\em Doc.\ Math.}, 6:127--183 (electronic), 2001.

\bibitem[CEM01]{ChaEchMey:01}
J\'{e}r\^{o}me Chabert, Siegfried Echterhoff, and Ralf Meyer.
\newblock Deux remarques sur l'application de {B}aum-{C}onnes.
\newblock {\em C.\ R.\ Acad.\ Sci.\ Paris}, 332, S\'{e}rie 1:607--610, 2001.

\bibitem[CEOO03]{CEO:03}
J{\'e}r{\^o}me Chabert, Siegfried Echterhoff, and Herv{\'e} Oyono-Oyono.
\newblock Shapiro's lemma for topological {$K$}-theory of groups.
\newblock {\em Comment. Math. Helv.}, 78(1):203--225, 2003.

\bibitem[DG83]{DuGil:83}
Maurice~J. Dupr\'{e} and Richard~M. Gillette.
\newblock {\em Banach bundles, Banach modules and automorphisms of
  {$C^*$}-algebras}.
\newblock Pitman Books Limited, 1983.

\bibitem[Dix64]{Dixmier:64}
Jacques Dixmier.
\newblock {\em Les {$C\sp{\ast} $}-alg\`ebres et leurs repr\'esentations}.
\newblock Cahiers Scientifiques, Fasc. XXIX. Gauthier-Villars \& Cie,
  \'Editeur-Imprimeur, Paris, 1964.

\bibitem[EKQR06]{EKQR:06}
Siegfried Echterhoff, S.~Kaliszewski, John Quigg, and Iain Raeburn.
\newblock A categorical approach to imprimitivity theorems for {$C\sp
  *$}-dynamical systems.
\newblock {\em Mem. Amer. Math. Soc.}, 180(850):viii+169, 2006.

\bibitem[FD88]{FellDo:88:1}
James M.~G. Fell and Robert~S. Doran.
\newblock {\em Representations of {$\sp *$}-algebras, locally compact groups,
  and {B}anach {$\sp *$}-algebraic bundles. {V}ol. 1}, volume 125 of {\em Pure
  and Applied Mathematics}.
\newblock Academic Press Inc., Boston, MA, 1988.
\newblock Basic representation theory of groups and algebras.

\bibitem[JT91]{JenThom:91}
Kjeld~Knudsen Jensen and Klaus Thomsen.
\newblock {\em Elements of {$KK$}-Theory}.
\newblock Birkh\"{a}user, 1991.

\bibitem[Laf02]{Lafforgue:02}
Vincent Lafforgue.
\newblock {$K$}-th\'{e}orie bivariante pour les alg\`{e}bres de {B}anach et
  conjecture de {B}aum-{C}onnes.
\newblock {\em Invent.\ Math.}, 149:1--95, 2002.

\bibitem[Laf06]{Lafforgue:06}
Vincent Lafforgue.
\newblock {$K$}-th\'{e}orie bivariante pour les alg\`{e}bres de {B}anach,
  groupo\"{i}des et conjecture de {B}aum-{C}onnes. {A}vec un appendice
  d'{H}erv\'{e} {O}yono-{O}yono.
\newblock {\em J.\ Inst.\ Math.\ Jussieu}, 2006.
\newblock Published online by Cambridge University Press 28 Nov 2006.

\bibitem[MRW87]{MuhReWill:87}
Paul~S. Muhly, Jean~N. Renault, and Dana~P. Williams.
\newblock Equivalence and isomorphism for groupoid {$C\sp \ast$}-algebras.
\newblock {\em J. Operator Theory}, 17(1):3--22, 1987.

\bibitem[OO98]{OyonoOyono:97}
Herv{\'e} Oyono-Oyono.
\newblock {\em La conjecture de {B}aum-{C}onnes pour les groupes agissant sur
  les arbres}.
\newblock PhD thesis, Universit\'{e} Claude Bernard - Lyon 1, 1998.

\bibitem[Par07]{Paravicini:07}
Walther Paravicini.
\newblock {\em $\KK$-Theory for {B}anach Algebras And Proper Groupoids}.
\newblock PhD thesis, Universit\"{a}t M\"{u}nster, 2007.
\newblock Persistent identifier: urn:nbn:de:hbz:6-39599660289.

\bibitem[Par08]{Paravicini:07:CNullX:erschienen}
Walther Paravicini.
\newblock A note on {B}anach {$\Cont_0(X)$}-modules.
\newblock {\em M{\"u}nster J. of Math.}, 1:267--278, 2008.
\newblock Persistent identifier: urn:nbn:de:hbz:6-43529451393.

\bibitem[Par09]{Paravicini:07:Induction:arxiv}
Walther Paravicini.
\newblock Induction for {B}anach algebras, groupoids and {KK}$^{\text{ban}}$.
\newblock {\em http://arxiv.org/abs/0902.4199}, 2009.
\newblock Preprint.

\bibitem[Pim86]{Pimsner:86}
Mihai~V. Pimsner.
\newblock {$KK$}-groups of crossed products by groups acting on trees.
\newblock {\em Invent. Math.}, 86(3):603--634, 1986.

\bibitem[Ren80]{Renault:80}
Jean~N. Renault.
\newblock {\em A Groupoid Approach to {C$^*$}-Algebras}, volume 793.
\newblock Springer-Verlag, Berlin, 1980.
\newblock Lecture Notes in Mathematics.

\end{thebibliography}
%

\end{document}